# On Isomorphism of Graphs and the *k*-clique Problem

Dhananjay P. Mehendale
Sir Parashurambhau College, Tilak Road, Pune-411030,
India

## Abstract

In this paper we develop three characterizations for isomorphism of graphs. The first characterization is obtained by associating certain bitableaux with the graphs. We order these bitableaux by suitably defined lexicographic order and denote the bitableau that is least in this order as the standard representation for the associated graph. The standard representation characterizes graphs uniquely. The second characterization is obtained in terms of associated rooted, unordered, pseudo trees. We show that the isomorphism of two given graphs is implied by the isomorphism of their associated pseudo trees. The third characterization is obtained in terms of ordered adjacency lists to be associated with two given labeled graphs. We show the two given labeled graphs are isomorphic if and only if their associated ordered adjacency lists can be made identical by applying suitable transpositions on any one of these lists. We discuss in brief the complexity of these characterizations described in this paper for deciding isomorphism of graphs. Finally, we discuss the *k*-clique problem in the light of these characterizations towards the end of the paper.

1. **Introduction:** When there exists an adjacency preserving bijection between their vertex sets the given graphs *G* and *H* are isomorphic. To determine whether two given graphs are isomorphic is called the Graph Isomorphism Problem (GI). GI is of great interest to computer scientists and researchers in other fields such as chemistry, switching theory, information retrieval, and linguistics. In particular GI is of profound interest to complexity theorists because yet the graph isomorphism problem is neither proved P nor proved NP-complete.
   The *k*-clique problem asks whether a *k*-clique can be found as a subgraph in the given graph *G*. It has been shown in the literature that the *k*-clique problem is NP-complete [1].
   Using the first isomorphism criterion, we begin with showing that testing isomorphism of two given graphs is equivalent to checking that their associated standard representations are identical. We propose a problem for further investigation of devising a smart procedure to find standard representation for graphs. We further see
using second criterion for isomorphism that solving the isomorphism
problem for two given graphs is in fact P. It follows from its
equivalence to the problem of isomorphism for two associated pseudo
trees. We show that problem of pseudo tree isomorphism checking has
same complexity as that of the problem of tree isomorphism checking.
The main idea behind first criterion can be summarized in terms of the
following first among the common notions of Euclid [2]: "*Things
which are equal to the same thing are also equal to one another*"

2. **Preliminaries:** Let $G$ be a $(p, q)$ graph, i.e. a graph on p points
   (vertices) and q lines (edges) with the following vertex set $V(G)$ and
   edge set $E(G)$ respectively:
   $$V(G) = \{v_1, v_2, \cdots, v_p\} \text{ and }$$
   $$E(G) = \{e_1, e_2, \cdots, e_q\}$$

**Definition 2.1:** The vertex adjacency bitableau associated with every
labeled copy of graph $G$, $VAB(G)$, is the following bitableau:

$$VAB(G) = \begin{pmatrix} 1 & \alpha_1^1 & \alpha_2^1 & \cdots & \alpha_{l_1}^1 \\ 2 & \alpha_1^2 & \alpha_2^2 & \cdots & \alpha_{l_2}^2 \\ \vdots & & & & \\ j & \alpha_1^j & \alpha_2^j & \cdots & \alpha_{l_j}^j \\ \vdots & & & & \\ p & \alpha_1^p & \alpha_2^p & \cdots & \alpha_{l_p}^p \end{pmatrix}$$

where left tableau represents the suffixes of the vertex labels and stand
for the vertices while the right tableau represents the rows of the
suffixes of the vertex labels and represent the vertices that are adjacent
to vertex whose suffix is written in the same row in the left tableau,
i.e. the appearance of entry $\alpha_k^j$ in the $j^{th}$ row of the right tableau
implies that vertex $v_j$ is adjacent to vertex $v_{\alpha_k^j}$.

For a given unlabeled $(p, q)$ graph $G$ there are $\dfrac{p!}{|\Gamma(G)|}$
distinct labeled copies of $G$, where $|\Gamma(G)|$ is the cardinality of the
automorphism group of $G$, [3]. With the right tableau of the $VAB(G)$

corresponding to each labeled copy of *G* we associate a sequence of numbers Seq ([m, n]), m≤ *p*, n ≤ *p* and thus order these labeled copies by the so called [m, n]-order. Thus,

Seq ([m, n]) = [1, 1] [2, 1] [3, 1] ... [*p*, 1] [1, 2] [2, 2] ... [*p*, 2] ... [*p*, *p*], where [m, n] = cardinality of number of entries ≤ m in the first n rows of the corresponding right tableau of *VAB*(*G*) of the labeled copy under consideration. It is easy to see that [m, n]-order is **a total order**. So, there exists a unique (up to automorphism) labeled copy of G for which the associated sequence has largest order.

We call the labeled copy of *G*, for which the associated sequence has largest [m, n]-order, the standard representation of *G* and the corresponding vertex adjacency bitableau the standard vertex adjacency bitableau, and denote it symbolically as *Std*(*VAB*(*G*)). It is clear that *VAB*(*G*) is an alternative representation of the adjacency matrix. In similar way we can associate the so called incidence bitableau, *IB*(*G*), with a graph which will stand for an alternative representation of the incidence matrix. Thus,

**Definition 2.2:** The incidence bitableau for graph *G*, *IB*(*G*), is the following bitableau:

$$IB(G) = \begin{pmatrix} 1 & \alpha_1^1 & \alpha_2^1 & \cdots & \alpha_{l_1}^1 \\ 2 & \alpha_1^2 & \alpha_2^2 & \cdots & \alpha_{l_2}^2 \\ \vdots & & & & \\ j & \alpha_1^j & \alpha_2^j & \cdots & \alpha_{l_j}^j \\ \vdots & & & & \\ p & \alpha_1^p & \alpha_2^p & \cdots & \alpha_{l_p}^p \end{pmatrix}$$

where left tableau represents the suffixes of the vertex labels and so represent the vertices while the right tableau represents the rows of the suffixes of the edge labels and so representing the edges that are incident on vertex whose suffix is written in the same row in the left tableau, i.e. the appearance of entry $\alpha_k^j$ in the *j*<sup>th</sup> row of the right tableau implies that vertex $v_j$ is incident on edge $e_{\alpha_k^j}$.

As is done previously, with the right tableau of the *IB*(*G*) corresponding to each labeled copy of *G* we associate a sequence of numbers

Seq ([m, n]), m≤ q, n ≤ p and thus order these labeled copies by the so called [m, n]-order. Thus,

Seq ([m, n]) = [1, 1] [2, 1] [3, 1] ... [q, 1] [1, 2] [2, 2] ... [q, 2] ... [q, p], where [m, n] = cardinality of number of entries ≤ m in the first n rows of the corresponding right tableau of the labeled copy under consideration. Since [m, n]-order is a total order, so, there exists a unique (up to automorphism) labeled copy of G for which the associated sequence has largest order.

We call the labeled copy of *G*, for which the associated sequence has largest order, the standard representation of *G* and the corresponding incidence bitableau the standard incidence bitableau, and denote it as *Std (IB(G))*.

**Algorithms for Standardization:** The above theorems will be of no practical use unless supplemented by some smart procedure to obtain the standard bitableaux, either *Std(VAB(G))* or *Std(IB(G))*, directly and efficiently, from a given unlabeled or labeled copy of *G*. We proceed to accomplish this task in this section.

In the language of matrices when two graphs *G* and *H* are isomorphic the corresponding isomorphism map induces a permutation $\sigma$ with permutation matrix M such that $A(G) = M A(H) M^T$, where A(*G*) and A(*H*) are the adjacency matrices for the graphs *G* and *H* respectively and $M^T$ is transpose of M. We now proceed to discuss

**The Action of the Permutation $\sigma$ on *VAB(G)*:** Since every permutation can be expressed as a product of transpositions or uniquely as the product of disjoint cycles it is enough to understand the action of transpositions or cycles on *VAB(G)*.

**Definition 2.3:** Let $(i, j)$ be a transposition, $i \neq j$ and $i \leq p, j \leq p$ then the action of this transposition on *VAB(G)*, [(i, j)] *VAB(G)*, involves the following operations:
(a) Interchanging of the rows *i* and *j* in the right tableau of *VAB(G)*.
(b) Replacing of the entries *i* by *j* and vice versa everywhere in the right tableau of *VAB(G)*.

**Definition 2.4:** Let $(i_1, i_2, \cdots, i_r)$ be a cycle (permutation) then the action of this cycle on *VAB(G)*, $[(i_1, i_2, \cdots, i_r)]$ *VAB(G)*, involves the following operations:

(a) Replacing of the rows $i_1$ by $i_2$, $i_2$ by $i_3$, ..., $i_r$ by $i_1$ in the right tableau of $VAB(G)$.

(b) Replacing of the entries $i_1$ by $i_2$, $i_2$ by $i_3$, ..., $i_r$ by $i_1$ everywhere in the right tableau of $VAB(G)$.

**Algorithm 2.1(Standardization of *VAB(G)*):**

(1) Label the unlabeled graph or alter by suitable transpositions the labels of the labeled graph $G$ such that the degrees of it's vertices $v_i$, say $d(v_i)$, satisfy $d(v_1) \geq d(v_2) \geq \cdots d(v_p)$ and form the corresponding $VAB(G)$.

(2) Apply the suitable transpositions $[(i, j)]$, with $d(v_i) = d(v_j)$, so that the nonincreasing order of the degrees is maintained and the first row of thus transformed $VAB(G)$ will contain smallest possible numbers maintaining degree sequence achieved in the first step.

(3) Now, apply the suitable transpositions so that the first row of the transformed $VAB(G)$ remains invariant and smallest possible numbers, maintaining the invariance of the first row and maintaining the degree sequence fixed in the first step, appear in the second row. (Note that if the first row after the execution of step (2) becomes $1 \mid j_1 j_2 \cdots j_m \cdots j_n \cdots j_{d(v_1)}$ then the first row remains invariant under transpositions $[(j_m, j_n)]$, such that $1 \leq m, n \leq d(v_1)$, and transpositions $[(k, l)]$ such that $k, l \neq j_p$, $1 \leq p \leq d(v_1)$.)

(4) Continue applying the transpositions which belong to the intersection of the set of invariant transpositions arrived at after achieving the smallest possible numbers in the earlier rows of the transformed $VAB(G)$, till we reach the last row where the process stops giving rise to $Std(VAB(G))$.

□

**Remark 2.1:** During the isomorphism check for the two given graphs by application of the standardization procedure if their vertex adjacency bitableaux become identical before we reach the corresponding standard bitableaux then we can stop declaring them isomorphic if we are only interested in the isomorphism question.

**Remark 2.2:** The necessary part of the necessary-sufficient condition proposed in the theorem 2.1 above can be effectively used to decide the non-isomorphism of two arbitrary graphs under consideration. When two arbitrary graphs are given for checking whether they are isomorphic, the usual approach is to try to show if they are not isomorphic [4]. This is done by checking whether some necessary condition is violated by these graphs. For checking the non-isomorphism of two arbitrary graphs using the necessary part of the necessary-sufficient condition provided by the theorem 2.1, it is not necessary to carry out complete standardization. We can very well stop the process of standardization started from the first row of the right tableaux of the corresponding vertex adjacency bitableaux for the graphs (algorithm 2.1), at the first occurrence of distinctness in the so far standardized rows from top row downwards and declare that the graphs are non-isomorphic.

**Example:** Consider following graphs, [5].

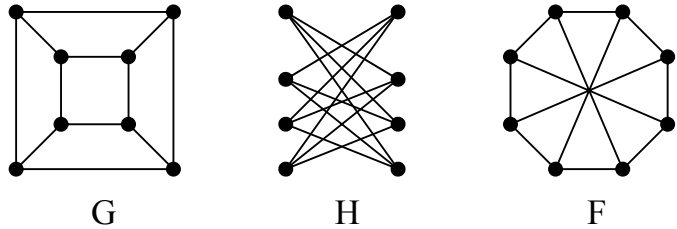

G　　　　　H　　　　　F

By application of standardization procedure we can find standard representations for graphs G, H, and F and can check that standard representations for graphs G, H agree while standard representation for F differs. So, G ≅ H but F is a different type of graph not isomorphic to *G* or H.

**The Action of a Permutation Pair ($\sigma, \tau$) on *IB(G)*:** Let G be a (*p, q*) graph. Note that the left tableau of *IB(G)* consists of *p* vertex labels while the right tableau consists of the rows of the edge labels incident on the vertices whose labels are written in the left tableau in the same rows. Also, vertex and edge labeling is independent and so can be assigned or altered independently. So, we take pair ($\sigma, \tau$) such that $\sigma \in S_p$ the group of permutations on *p* symbols while $\tau \in S_q$ the group of permutations on *q* symbols, and apply them respectively on the left and the right tableau of *IB(G)*.

　　　　　We now proceed to define the action of permutations on *IB(G)*.

**Definition 2.5:** Let $(i, j)$ be a transposition, $i \neq j$ and $i \leq p, j \leq p$ then the action of this transposition on the left tableau of $IB(G)$, $[(i,j)]IB(G)$, involves the following operation:

(a) Interchanging of the rows $i$ and $j$ in the right tableau of $IB(G)$.

**Definition 2.6:** Let $(i_1, i_2, \cdots, i_r)$ be a cycle (permutation) then the action of this cycle on the left tableau of $IB(G)$, $[(i_1, i_2, \cdots, i_r)]IB(G)$, involves the following operations:
(a) Replacing of the rows $i_1$ by $i_2$, $i_2$ by $i_3$, ..., $i_r$ by $i_1$ in the right tableau of $IB(G)$.

**Definition 2.7:** Let $(i, j)$ be a transposition, $i \neq j$ and $i \leq p, j \leq p$ then the action of this transposition on the **right** tableau of $IB(G)$, $[(i,j)]IB(G)$, involves the following operation:
(a) Replacing of the entries $i$ by $j$ and vice versa everywhere in the right tableau of $IB(G)$.

**Definition 2.8:** Let $(i_1, i_2, \cdots, i_r)$ be a cycle (permutation) then the action of this cycle on the right tableau of $IB(G)$, $[(i_1, i_2, \cdots, i_r)]IB(G)$, involves the following operations:
(a) Replacing of the entries $i_1$ by $i_2$, $i_2$ by $i_3$, ..., $i_r$ by $i_1$ everywhere in the right tableau of $IB(G)$.

Hereafter, when we say that the action of transposition (permutation) is done on left tableau then it is the action on the vertex labels and when we say that the action of transposition (permutation) is done on right tableau then it is the action on the edge labels.

**Algorithm 2.2(Standardization of $IB(G)$):**

(1) Label the unlabeled graph or alter by suitable transpositions the vertex labels of the labeled graph $G$ such that the degrees of its vertices $v_i$, say $d(v_i)$, satisfy $d(v_1) \geq d(v_2) \geq \cdots d(v_p)$. Also,
(2) Apply the suitable transpositions on edge labels in the right tableau of $IB(G)$, if and when required, so that the first row of thus transformed $IB(G)$ becomes $1| 1\ 2\ 3\ ...j\ ...\ d(v_1)$.

(3) Now, apply the suitable transpositions on right tableau so that the first row of the transformed *IB(G)* remains invariant and smallest possible numbers, maintaining the invariance of the first row and maintaining the degree sequence fixed in the first step, appear in the second row.

(4) Continue applying the transpositions which belong to the intersection of the set of invariant transpositions arrived at after achieving the smallest possible numbers in the earlier rows of the transformed *IB(G)*, till we reach the last row where the process stops giving rise to *Std(IB(G))*.

□

**Remark 2.3:** It is easy to check that it is possible to achieve step (2) and when step (2) is accomplished this row will remain invariant under the further action of transpositions $(i, j)$ on the right tableau, $i \neq j$ when both $1 \leq i \leq d(v_1)$, $1 \leq j \leq d(v_1)$ or both $i > d(v_1)$, $j > d(v_1)$ hold.

### 3. Characterization in terms of *VAB(G)* and *IB(G)*:

**Theorem 3.1:** Let $G$ and $H$ be two $(p, q)$ graphs. $G \cong H$ if and only if $Std(VAB(G)) = Std(VAB(H))$.

**Proof:** Follows from the uniqueness of the standard representation.

□

**Theorem 3.2:** Let $G$ and $H$ be two $(p, q)$ graphs. $G \cong H$ if and only if $Std(IB(G)) = Std(IB(H))$.

**Proof:** Follows from the uniqueness of the standard representation.

□

Thus, we have reduced the problem of checking the isomorphism of two given graphs to obtaining their standard bitableaux (either $Std(VAB(G))$ or $Std(IB(G))$) and check their equality.

### 4. Graph Isomorphism:
We obtain a polynomial time algorithm for testing isomorphism of two graphs in terms of isomorphism of the associated rooted, unordered, pseudo trees. It is well known that there exists a polynomial time algorithm for testing isomorphism of trees.

We see that such algorithm can be extended to check isomorphism of rooted, unordered, pseudo trees maintaining its polynomial time complexity.

**Pseudo Tree with Connected Graph:** We introduce now the idea of associating a rooted, unordered, pseudo tree with connected graph. We show here the use of this association to obtain a new polynomial time algorithm for testing isomorphism of graphs. Thus we show that isomorphism problem for graphs can be looked upon as isomorphism problem for its associated rooted, unordered, pseudo trees. We further show that the required algorithm for testing isomorphism of rooted, unordered, pseudo trees doesn't differ in complexity from the algorithm for testing isomorphism of trees (which can be done in quadratic time). We associate with (or, look at a) connected graph (as) a rooted, unordered, pseudo tree. To describe this idea let us begin with an

**Example:** Consider following connected graph, G say, and its associated rooted, unordered, pseudo tree, T say:

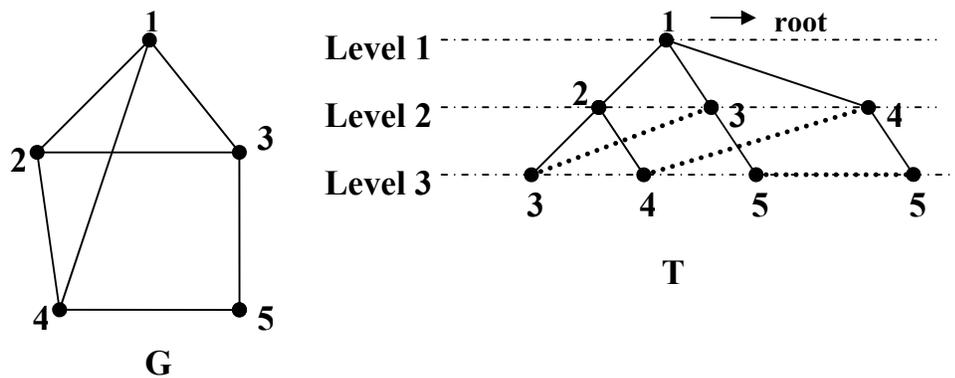

This tree T is called rooted because it has a root (vertex with label 1), unordered because it contains multiple vertices with same label, and pseudo because it contains pseudo edges joining vertices with same label shown by dotted line segments which are actually nonexistent.

**Algorithm 4.1 (Pseudo tree Construction Algorithm):** For a given unlabeled ($n$, $e$) graph we carry out following steps to construct rooted, unordered, pseudo tree.
1) Choose any vertex of given unlabeled graph as root and assign label "1" to it.
2) Assign labels "2", "3", …, "n" to other vertices of graph.
3) Create 1$^{st}$ level of desired tree consisting of vertex with label "1".
4) Take all vertices adjacent to vertex with label "1" and paint the edges emerging from vertex with label "1" and entering in all vertices s1, s2, … etc adjacent to it with **red** color in the graph as a distinguishing mark. This marking is done to understand that these edges are now used and are not to be used again.
5) Create 2$^{nd}$ level by putting all vertices "s1", "s2", … etc., adjacent to the only vertex with label "1" in 1st level, in the 2$^{nd}$ level and draw the corresponding edges joining vertex labeled "1" in the first level to each vertex, among the vertices "s1", "s2", … etc, in the second level.
6) Take vertex "s1" in the graph. Treat all vertices connected to it by a red edge as nonadjacent (e.g. in this case the vertex with label "1"). Determine all vertices adjacent to vertex "s1", say "t11", "t12", …, etc. and paint the edges emerging from vertex with label "s1" and entering in all vertices "t11", "t12", … etc adjacent to it with red color in the graph as a distinguishing mark. . This marking is done to understand that these edges are now used and are not to be used again.
7) Create 3$^{rd}$ level partially by putting vertices "t11", "t12", … etc., adjacent to the vertex with label "s1" in 2nd level, in the 3$^{rd}$ level and draw the corresponding edges joining vertex labeled "s1" in the second level to each vertex, among the vertices "t11", "t12", … etc, in the third level.
8) Repeat steps 6), 7) for all other vertices "s2", "s3", … etc. in the 2$^{nd}$ level, by taking them one by one in a sequence, and thus complete the creation of 3$^{rd}$ level, i.e. create the third level completely.
9) Continue this procedure of pseudo tree construction till all edges in the graph become red and all the three levels have been created completely.
10) In the rooted, unordered, pseudo tree since we avoid cycle formation by following special procedure for its construction described in above steps (which forcefully forbids presence of edges joining vertices in the same level), many vertices show their more than once appearance. So, finally, join all pairs of vertices with same label by a pseudo edge (dotted line segment).

□

Important care to be taken while we construct rooted, unordered, pseudo tree:

"*All vertices adjacent to root ($1^{st}$ level) must appear in the $2^{nd}$ level. All vertices adjacent to vertices in the $2^{nd}$ level connecting to so far unused and so unmarked edges must appear in the $3^{rd}$ level and so on and every circuit formation is avoided and finally all vertices with same label are joined to each other by dotted line segments representing pseudo edges which are actually nonexistent*"

**Remark 4.1:** Note that in order to avoid formation of any cycle (which should not exist in a tree) we take such vertices as additional vertices in the next level so clearly the rooted tree contains more vertices than the one contained in the original graph.

**Remark 4.2:** It is straightforward to check that if the connected graph G under consideration is an (n, e) graph, i.e. it contains n vertices and e edges, then its associated rooted, unordered, pseudo tree contains (e+1) vertices.

**Remark 4.3:** It is easy to see further that in all there are (e+1-n) repetitions of vertices, i.e. there are in all (e+1-n) vertices which show appearance more than once in the associated rooted, unordered, pseudo tree.

**Remark 4.4:** There are well known linear time algorithms to test isomorphism of rooted trees [6]. For problem of isomorphism in general case, i.e. unrooted trees, where we do not have a useful start as root, the algorithm for rooted trees can be easily adopted to make it work on unrooted trees: we only have to pick a vertex in first tree and declare it as root, and try all the vertices of second tree one after the other. The two unrooted trees will be isomorphic if and only if there is at least one vertex in second tree which, when declared as root prompts the linear time algorithm to test isomorphism of rooted trees to answer positively. This makes the algorithm for testing isomorphism of unrooted trees quadratic in time.

**Graph Isomorphism (Characterization in terms of Pseudo Trees):**
In this section we state the result that decides the isomorphism of two given graphs in $\sim O(n^2)$.

**Definition 4.1:** Two rooted, unordered, pseudo trees are isomorphic if they are isomorphic as rooted trees under some adjacency preserving bijection, $\sigma$ and further $\sigma$ also preserves pseudo-adjacency, i.e. if w is a vertex that occurs more than one time, say k times, in first pseudo tree so that there are pseudo edges joining each pair formed from these multiple copies of w, then the image of w, $\sigma(w)$, also will appear k times in second rooted, unordered, pseudo (image) tree and there are also pseudo edges present, joining each pair of these corresponding multiple copies of $\sigma(w)$, under the same bijection $\sigma$.

**Theorem 4.1:** Let G and H be the two given (n, e) graphs. Let V(G) and V(H) be the vertex sets for G and H respectively. Let vertex u belongs to V(G) and vertex v belongs to V(H) and let Tu(G) and Tv(H) be the associated rooted, unordered, pseudo trees rooted at vertex u and v for graphs G and H respectively. If Tu(G) is isomorphic to Tv(H) then the graph G is isomorphic to graph H under the same mapping depicting isomorphism of associated pseudo trees.

**Proof:** As we have Tu(G) is isomorphic to Tv(H) we have got adjacency and pseudo-adjacency preserving bijection $\sigma$. Now, if all the pseudo edges are contracted then we see that the number of vertices change from (e+1) to n and the same bijection, $\sigma$, will now act as isomorphism for graphs G and H from which these (isomorphic) pseudo trees were constructed.

□

We now proceed to discuss the algorithm for checking isomorphism of given two graphs which has polynomial order because the polynomial order is enough for checking isomorphism of rooted, unordered, pseudo trees associated with these graphs.

**Algorithm 4.2 (Graph Isomorphism Checking)):**

1) Let G and H be two (n, e) graphs given for isomorphism check. Let V(G) and V(H) be their vertex sets respectively. Let these sets are V(G) = {u1, u2, u3, ….} and V(H) = {v1, v2, v3, ….}.
2) Construct rooted, unordered, pseudo tree, say Tu1(G), rooted at vertex u1.

3) Construct rooted, unordered, pseudo tree, say Tv1(H), rooted at vertex v1.
4) Check whether Tu1(G) is isomorphic to Tv1(H).
5) If yes, declare that G is isomorphic to H and stop.
6) Else, take next vertex v2 in V(H) as root and construct rooted, unordered, pseudo tree, Tv2(H), rooted at vertex v2.
7) Check whether Tu1(G) is isomorphic to Tv2(H).
8) If yes, declare that G is isomorphic to H and stop.
9) Else, continue the above steps for next vertex v3, v4, …, till we reach either at isomorphism decision for G and H or at the last vertex of V(H).
10) If we don't find isomorphism of Tu1(G) with any of the rooted, unordered, pseudo trees, that we construct using vertices of V(H) by taking them one by one in succession as roots, then declare that G and H are not isomorphic.

□

**Example:** We consider below three graphs and find their associated pseudo trees. The isomorphism of these three graphs becomes evident from the easy isomorphism of their associated pseudo trees (pseudo edges and levels are not shown to avoid clumsiness).

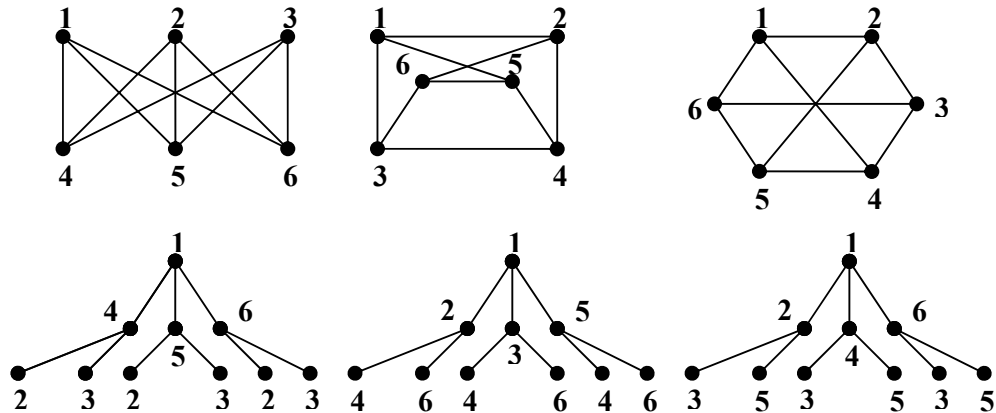

**Remark 4.5:** Since isomorphism of rooted, unordered, pseudo trees can be checked in linear time using any known algorithm and we may require to repeat this checking $n$ number of times ($n = |V(H)|$, cardinality of vertex set) for desired decision, we get the overall complexity for algorithm for isomorphism of graphs as $\sim O(n^2)$.

**Graph Isomorphism (Characterization using transpositions and reordering of ordered adjacency list):** We now proceed to discuss one more algorithm for Graph Isomorphism which in the worst case is of order $\sim O(e^3)$, where $e$ stands for number of edges in the graph. In practice this algorithm works very efficiently and actually doesn't require the worst case complexity, namely, $\sim O(e^3)$ to produce the decision about the isomorphism of given two ($n$, $e$) graphs $G$ and $H$. In this algorithm we may need to apply in the worst case the action of $e$ transpositions, and each such action of transposition further requires ordering of e edge labels using any standard algorithm of order $\sim O(e^2)$, thus, total order of our algorithm becomes $\sim O(e^3)$ in the worst case.

Let $G$ be a labeled ($p$, $q$) graph whose vertices are labeled with labels $\{1, 2, 3, \ldots, p\}$. If there is an edge joining two vertex labels $i$ and $j$ say such that $i < j$ then we associate pair $(i, j)$ as label with that edge.

**Definition 4.1:** The **ordered adjacency list** is the ordered list of edge labels presented in a row (or column) as follows: If vertex with label 1 is adjacent to vertices with labels $\{i_1, i_2, \cdots, i_{k_1}\}$ and $\{1 < i_1 < i_2 < \cdots < i_{k_1}\}$ then the first $k_1$ entries of adjacency list are $(1, i_1)(1, i_2)(1, i_3) \cdots (1, i_{k_1})$. If vertex with label 2 is adjacent to vertices with labels $\{j_1, j_2, \cdots, j_{k_2}\}$ and $\{2 < j_1 < j_2 < \cdots < j_{k_2}\}$ then the next $k_2$ entries of adjacency list are $(2, j_1)(2, j_2)(2, j_3) \cdots (2, j_{k_2})$. We then take in succession the vertices with the vertex labels 3, 4, …. and continue extending the list along same lines by the appending the next proper ordered pairs to the already constructed list till we reach the vertex with label ($p$-1) and finally by appending (or not appending) the vertex pair ($p$-1, $p$) when the vertex with label ($p$-1) is adjacent (nonadjacent) to vertex with label $p$ we thus complete construction of ordered adjacency list for the graph.

It is well known to all that two graphs are isomorphic if and only if there exists one-one adjacency preserving map between their vertex sets, i.e. one graph can be transformed into other by relabeling the vertices of (any) one graph. More formally, let $G$ and $H$ be two ($p$, $q$) graphs and let $A_G$, $A_H$ be the adjacency matrices associated with graphs $G$ and $H$ respectively then these graphs $G$ and $H$ are isomorphic if and only if there exists a permutation matrix $P$ such that $A_G = PA_H P^{-1}$.

Now the problem with using this result to test two given graphs for isomorphism is actually as follows: If given graphs are isomorphic how to discover such permutation $\sigma$ whose associated permutation matrix, $P$ say, satisfies $A_G = PA_H P^{-1}$? Also, when given graphs will not be isomorphic there will not exists any such permutation $\sigma$ whose associated permutation matrix, $P$ say, achieves $A_G = PA_H P^{-1}$. How to conclusively arrive at the result that such permutation won't exist when two given graphs under consideration are nonisomorphic?

The answer to these questions can be extracted from the simple fact that any permutation (and so also the desired permutation) is product (composition) of transpositions. We will see that the desired permutation builds up from composition of suitable transpositions to be applied in succession on the ordered adjacency list of any one of the two given graphs. When the two given graphs are isomorphic one can write down their associated ordered adjacency lists. One then choose (any) one list and choose proper transpositions to be applied on this list for making it identical with the other ordered adjacency list. When the two given graphs are not isomorphic then one cannot make the lists identical by whatever transpositions one will choose to apply on (any) one ordered adjacency list to make it identical with the other ordered adjacency list, i.e. when one chooses the next transposition for further identification of two ordered adjacency lists, which have been made partially identical by earlier transpositions, then the new transposition disturbs the already obtained partial identification of the ordered adjacency lists.

**Action of transposition on given ordered adjacency list:** Let $L$ be an ordered adjacency list associated with some labeled ($p$, $q$) graph $G$ and its vertices are labeled with numbers $\{1, 2, 3, ..., p\}$. Let us denote the transposition by $[i, j]$. The action of this transposition on ordered adjacency list, expressed as $[i, j](L)$, is a two step procedure:
(i) It changes everywhere in the ordered adjacency list $L$ the symbol $i$ by symbol $j$ and vice versa.

**(ii)** After the above step it reorders the (modified) adjacency list so that again it becomes ordered in the sense of definition 4.1.

**Theorem 4.2:** Let $G$ and $H$ be two $(p, q)$ graphs. $G \cong H$ if and only if there exists a permutation made up of composition of transpositions which when applied successively on ordered adjacency list of $G$ say, converts it successfully into ordered adjacency list of $H$.

**Proof:** Suppose $G$ is isomorphic to $H$ then by definition there exists permutation $\sigma$ whose associated permutation matrix, $P$ say, satisfies $A_G = PA_H P^{-1}$. If we break the permutation $\sigma$ into transpositions and will apply these transpositions successively on ordered adjacency list of $H$ we will straightway get the ordered adjacency list of $G$.
Now, suppose $G$ is not isomorphic to $H$ then there will not exist permutation $\sigma$ whose associated permutation matrix, $P$ say, satisfies $A_G = PA_H P^{-1}$ so if any $\sigma$ will be taken and will be broken into transpositions and applied successively on ordered adjacency list of $H$ it will certainly fail to produce the ordered adjacency list of $G$.

□

**Algorithm 4.3 (Isomorphism using ordered adjacency list):**

1) Take labeled copies of given $(p, q)$ graphs $G$ and $H$ labeled with labels $\{1, 2, \ldots, p\}$ and prepare ordered adjacency lists say $L_G$ and $L_H$ and select $L_H$ for applying transpositions and to see whether we can equalize $L_H$ ultimately with $L_G$ by action of successive suitable transpositions.

2) Let the first elements in the ordered adjacency lists $L_G$ and $L_H$ be $(1, i_1)$ and $(1, j_1)$ respectively. If $i_1 = j_1$ then proceed to next elements in the ordered adjacency lists $L_G$ and $L_H$, namely, $(1, i_2)$ and $(1, j_2)$ to check whether $i_2 = j_2$. Else, if $i_1 \neq j_1$ then carry out action of transposition on $L_H$ that replaces $i_1$ by $j_1$ and vice versa everywhere in $L_H$ and then order the changed $L_H$ in the sense of definition 4.1. This new $L_H$ now will stand for $L_H$ for further actions of transpositions. We represent this (two stepped) action, carried out when $i_1 \neq j_1$, symbolically as $[i_1, j_1](L_H)$. Note that this step has equalized the first element in the ordered adjacency lists $L_G$ and (new) $L_H$ which will now stand for $L_H$.

3) Now, proceed to next elements in the ordered adjacency lists $L_G$ and (new) $L_H$, namely, $(1, i_2)$ and $(1, j_2)$ and check whether $i_2 = j_2$. If $i_2 = j_2$ then proceed to next elements in the ordered adjacency lists $L_G$ and $L_H$, namely, $(1, i_3)$ and $(1, j_3)$ to check whether $i_3 = j_3$.. Else, if $i_2 \neq j_2$ then carry out action of transposition on $L_H$ that replaces $i_2$ by $j_2$ and vice versa everywhere in $L_H$ and then as previous order the changed $L_H$ in the sense of definition 4.1. This new $L_H$ now will stand for $L_H$ for further actions of transpositions. We represent this (two stepped) action that carried out symbolically as $[i_2, j_2](L_H)$. Note that if action $[i_2, j_2](L_H)$ performed to equalize second elements in the ordered adjacency lists $L_G$ and $L_H$ does not disturbs the already achieved equality of first elements in the lists then up to this step one has equalized the first two elements in the ordered adjacency lists $L_G$ and (new) $L_H$ which will now stand for $L_H$. .

4) We continue in this way with equalizing next element in $L_H$ with the element in $L_G$ in the same place by the action of suitably chosen transpositions on (new) $L_H$ which now stands for $L_H$ till we can proceed along these lines without disturbing the equality of elements achieved earlier.

$\square$

**Remark 4.6:** Note that when given $(p, q)$ graphs $G$ and $H$ are isomorphic we should be able to achieve equality of $L_G$ and $L_H$ and when we are unable to achieve equality of $L_G$ and $L_H$ by algorithm 4.3 then $G$ and $H$ are not isomorphic.

**Example:** We now proceed to discuss one example to apply this algorithm. The graphs in this example are isomorphic. A similar example can be considered for the case of nonisomorphic graphs. There one will see that earlier obtained partial equality of ordered adjacency lists cannot be maintained in attempt of extending the equality further. Consider following labeled graphs *G, H, K* given below:

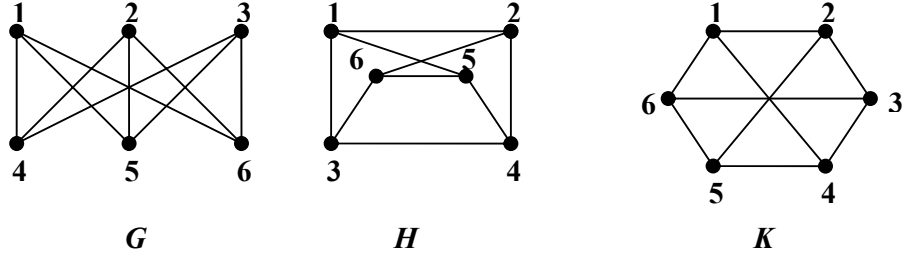

*G*  *H*  *K*

We now proceed to record the ordered adjacency lists $L_G, L_H$, and $L_K$. We then carry out action of suitable transpositions, firstly on $L_H$ and then on $L_K$, and show that we can identify $L_H, L_K$ with $L_G$. This shows that graphs *G, H, K* are isomorphic.

$$L_G = (1,4)(1,5)(1,6)(2,4)(2,5)(2,6)(3,4)(3,5)(3,6)$$

$$L_H = (1,2)(1,3)(1,5)(2,4)(2,6)(3,4)(3,6)(4,5)(5,6)$$

$$L_K = (1,2)(1,4)(1,6)(2,3)(2,5)(3,4)(3,6)(4,5)(5,6)$$

We first check whether $L_H$ can be made identical to $L_G$ by the action of properly chosen transpositions: Consider

$$[2,4](L_H) = (1,3)(1,4)(1,5)(2,3)(2,4)(2,5)(3,6)(4,6)(5,6) = L_H^1$$

Here we now set $L_H^1 = L_H$. Proceeding with this (new) $L_H$ we do

$$[3,6](L_H) = (1,4)(1,5)(1,6)(2,4)(2,5)(2,6)(3,4)(3,5)(3,6) = L_G$$

Let us now proceed to apply suitable transpositions on $L_K$. Consider

$$[2,5](L_K) = (1,4)(1,5)(1,6)(2,4)(2,5)(2,6)(3,4)(3,5)(3,6) = L_G$$

Therefore, graphs *G, H, K* are isomorphic.

5. **The *k*-clique Problem:** To deal with this problem we make the following change in the standardization procedure. It essentially consists of relaxing the condition of maintaining the nonincreasing order of the degree sequence and restricting it to degree (*k*−1) when we are seeking the existence of a subgraph isomorphic to a *k*-clique in the given graph. Also, the [m, n] order is restricted to degree (*k*−1). To make the things precise we state below few definitions:

**Definition 5.1:** The (*k*−1) degree restricted (restrictedness indicated by a vertical bar after (*k*−1) numbers in each row) vertex adjacency bitableau associated with every labeled copy of graph *G*, *VAB(G)*, is the following bitableau:

$$VAB(G, k-1) = \begin{pmatrix} 1 & \alpha_1^1 & \cdots & \alpha_{(k-1)}^1 & \cdots & \alpha_{l_1}^1 \\ 2 & \alpha_1^2 & \cdots & \alpha_{(k-1)}^2 & \cdots & \alpha_{l_2}^2 \\ \vdots & & & & & \\ j & \alpha_1^j & \cdots & \alpha_{(k-1)}^j & \cdots & \alpha_{l_j}^j \\ \vdots & & & & & \\ p & \alpha_1^p & \cdots & \alpha_{(k-1)}^p & \cdots & \alpha_{l_p}^p \end{pmatrix}$$

where (as in the case of *VAB(G)*) left tableau represents the suffixes of the vertex labels and so stand for the vertices while the right tableau represents the rows of the suffixes of the vertex labels and represent the vertices that are adjacent to vertex whose suffix is written in the same row in the left tableau, i.e. the appearance of entry $\alpha_k^j$ in the $j^{th}$ row of the right tableau implies that vertex $v_j$ is adjacent to vertex $v_{\alpha_k^j}$.

With the right tableau of the *VAB(G, k−1)* corresponding to each labeled copy of *G* we associate a sequence of numbers Seq ([m, n]), m ≤ p, n ≤ p and thus order these labeled copies by the so called (*k*−1) degree restricted [m, n]-order (by restricted we meant that the numbers only up to and including first (k−1) numbers in each row up to first *k* rows are considered and the other numbers are ignored in the counting). Thus,

Seq ([m, n]) = [1, 1] [2, 1] [3, 1] ... [*p*, 1] [1, 2] [2, 2] ... [*p*, 2] ... [*p*, *p*],

where [m, n] = cardinality of number of entries ≤ m counted only up to and including first (k−1) numbers of each row in the first k rows of the corresponding right tableau of $VAB(G, k-1)$ of the labeled copy under consideration. It is easy to see that this (k−1) degree restricted [m, n]-order is a total order. So, there exists a unique (up to automorphism of the subgraph represented by the leading k by $(k-1)$ sub-bitableau) labeled copy of G for which the associated sequence has largest order.

We call the labeled copy of G, for which the associated sequence has largest restricted [m, n]-order, the (k−1) degree restricted standard representation of G and the corresponding vertex adjacency bitableau the (k−1) degree restricted standard vertex adjacency bitableau, and denote it symbolically as $Std(VAB(G, k-1))$.

**Theorem 5.1:** Let G be the given (p, q) graph. G has a k-clique as a subgraph if and only if $Std(VAB(G, k-1))$ contains the following bitableu as a subbitableau at the leading position:

$$\begin{pmatrix} 1 & 2 & 3 & \cdots & & k \\ 2 & 1 & 3 & \cdots & & k \\ \vdots & & & & & \\ j & 1 & 2 & \cdots & & k \\ \vdots & & & & & \\ k & 1 & 2 & \cdots & & k-1 \end{pmatrix}$$

**Remark 5.1:** It is clear to see that the above subbitableau is nothing but the standard bitableau corresponding to a k-clique, and also every bitableau representing a complete graph is inherently standard and nothing is to be done for its standardization.

**Proof:** When the above bitableau will be present the existence of a k-clique as a subgraph is clear. When k-clique as a subgraph is present in the given graph it must show its appearance in the form of the above bitableau at the leading position due to uniqueness of the (k−1) degree restricted standard representation contained as a subbitableau.

□

**Algorithm 5.1(Restricted Standardization of *VAB(G, k−1)*):**

(2) Label the unlabeled graph $G$ with vertex labels $\{1, 2, \ldots, p\}$ and by suitable transpositions rename the labels of the labeled graph $G$ such that the degrees of it's vertices $v_i$, say $d(v_i)$, satisfy $d(v_1) \geq d(v_2) \geq \cdots d(v_p)$ and form the corresponding *VAB(G, k−1)* and apply suitable transpositions $[(i, j)]$, among the vertices of $G$ having degrees $\geq (k-1)$ so that the first row of thus transformed *VAB(G, k−1)* will contain smallest possible numbers up to first $(k-1)$ numbers.

(3) Now, apply the suitable transpositions which maintain the degrees of at least the first $k$ rows $\geq (k-1)$, keeps the numbers up to first $(k-1)$ numbers in the first row of the transformed *VAB(G, k−1)* invariant, and cause the appearance of smallest possible numbers in the second row, when numbers are considered up to first $(k-1)$ numbers, maintaining the invariance of the first row. (Note: if the first row after the execution of step (2) becomes $1 \mid j_1 j_2 \cdots j_m \cdots j_n \cdots j_{d(v_1)}$ then the first $(k-1)$ numbers in the first row remain invariant under transpositions $[(j_m, j_n)]$, such that $1 \leq m, n \leq (k-1)$, and transpositions $[(k, l)]$ such that $k, l \neq j_p$, $1 \leq p \leq (k-1)$.)

(4) Continue applying the transpositions among the vertices of $G$ having degrees $\geq (k-1)$ which belong to the intersection of the set of invariant transpositions (for the first $(k-1)$ numbers) arrived at after achieving the smallest possible numbers up to first $(k-1)$ numbers in each of the earlier rows of the transformed bitableau, till we reach the $k$-th row where the process stops giving rise to a above mentioned subbitableau corresponding to a $k$-clique if and when a $k$-clique exists as a subgraph of $G$.

□

**Remark 5.3:** It is easy to check that rooted, unordered, pseudo tree associated with a complete graph contains exactly three levels, i.e. it is always a three leveled typical tree. Because of total symmetry of complete graphs we can take any vertex as root and build the associated pseudo tree. For example for complete graph on n vertices with V(G) = $\{1, 2, 3, \ldots ,n\}$ the structure of this tree is as follows: The first level contains vertex with label 1 called root. The second level contains vertices $\{2, 3, \ldots , n\}$ joined to root by an edge. The third level separately contains vertices $\{3, 4, \ldots, n\}$ joining vertex with label 2 in second level, vertices $\{4, 5, \ldots, n\}$ joining vertex with label 3 in second level, ……, vertex {n} joining vertex with label (n-1) in second level. Further, there will be present the pseudo edges joining all the pairs of vertices with repeated occurrence.

**Remark 5.4:** The problem of checking the existence of *k*-clique thus reduces to checking the existence special kind of rooted, unordered, pseudo sub-tree corresponding to *k*-clique described in the above remark 5.3 having vertex with label 1 or some other vertex as root.

We now conclude with following problem for further investigations:

**Problem:** Devise efficient algorithms to obtain standard representation for various bitableaux defined in this paper.

### Acknowledgements

I am thankful to Prof. M. R. Modak for useful discussions.